\documentclass[12pt]{amsart}

\newtheorem{thm}{Theorem}
\newtheorem{prop}{Proposition}

\newtheorem*{theorem}{Theorem}
\newtheorem{cor}{Corollary}
\newtheorem{lem}{Lemma}
\newtheorem{conj}{Conjecture}
\newtheorem{rem}{Remark}

\newcommand{\N}{{\mathbb N}}

\newcommand{\ppp}{{\mathfrak P}}\newcommand{\pp}{{\mathcal P}}

\newcommand{\R}{{\mathbb R}}\newcommand{\RR}{{\mathbb R}^2}

\newcommand{\const}{\mbox{const}} 

\newcommand{\tx}{\tilde{x}}

\newcommand{\tM}{\tilde{M}}

\newcommand{\al}{\alpha}
\newcommand{\be}{\beta}
\newcommand{\ga}{\gamma}\newcommand{\Ga}{\Gamma}
\newcommand{\de}{\delta}
\newcommand{\ve}{\varepsilon}
\newcommand{\la}{\lambda}
\newcommand{\ka}{\kappa}
\newcommand{\vap}{\varphi}
\newcommand{\si}{\sigma}

\begin{document}

\bibliographystyle{plain}

\title[Counting and Blocking Geodesics]
{Topological entropy and blocking cost for geodesics in riemannian
manifolds}

\author{Eugene Gutkin}

\address{IMPA \\
Estrada Dona Castorina 110 \\
Rio de Janeiro,   22460-320\\
Brasil; UMK\\
Chopina 12/18\\
87 -- 100 Torun\\
and IMPAN\\
Sniadeckich 8\\
00-956 Warszawa\\
Poland}

\email{\ gutkin@mat.uni.torun.pl,\ gutkin@impa.br}



\date{\today}

\begin{abstract}
For a pair of points $x,y$ in a compact, riemannian manifold $M$
let $n_t(x,y)$ (resp. $s_t(x,y)$) be the number of geodesic
segments with length $\leq t$ joining these points (resp. the
minimal number of point obstacles needed to block them). We study
relationships between the growth rates of $n_t(x,y)$ and
$s_t(x,y)$ as $t\to\infty$. We derive lower bounds on $s_t(x,y)$
in terms of the topological entropy $h(M)$ and its fundamental
group. This strengthens the results of Burns-Gutkin~\cite{BG06}
and Lafont-Schmidt~\cite{LS}. For instance, by~\cite{BG06,LS},
$h(M)>0$ implies that $s$ is unbounded; we show that $s$ grows
exponentially, with the rate at least $h(M)/2$.
\end{abstract}

\maketitle

\section*{Introduction}       \label{intro}
By a riemannian manifold we will always mean a closed, complete,
connected, infinitely differentiable, riemannian manifold. Let $M$
be a riemannian manifold. By a geodesic $\gamma\subset M$ we will
mean an oriented geodesic segment; thus, $\gamma$ has endpoints
$x,y\in M$ and a positive length, $|\gamma|$. (We allow $x=y$.) If
$z\in M$ is an interior point of $\gamma$, we say that $\gamma$
{\em passes through $z$}. For $x,y\in M$ and $0<t$ we denote by
$G_t(x,y)$ the set of geodesics with endpoints $x,y$ and length at
most $t$. Let $\Gamma_t(x,y)\subset G_t(x,y)$ be the subset of
those $\gamma\in G_t(x,y)$ that do not pass through either $x$ or
$y$. We set
$G(x,y)=\cup_{t\in\R_+}G_t(x,y),\Ga(x,y)=\cup_{t\in\R_+}\Ga_t(x,y)$.
We will say that the geodesics in $G(x,y)$ (resp. $\Ga(x,y)$) {\em
join} (resp. {\em connect}) $x$ with $y$.

A finite set $B\subset M\setminus\{x,y\}$ is a {\em blocking set}
for $\Gamma_t(x,y)$ if every $\gamma\in\Gamma_t(x,y)$ passes
through a point in $B$. Let $s_t(x,y)\le\infty$ be the minimal
cardinality of a blocking set, and set $s(t)=\sup_{x,y\in
M}s_t(x,y)$. We say that $t\mapsto s_t(x,y)$ is the {\em blocking
threshold} function for $x,y\in M$, and that $t\mapsto s(t)$ is
the {\em blocking cost} function for $M$.

The framework of {\em security} for riemannian manifolds concerns
the question of blocking all of geodesics in $\Ga(x,y)$ by a
finite set \cite{BG06,GS06,PH,LS}. A pair $x,y\in M$ is secure if
there is a finite blocking set for $\Ga(x,y)$; otherwise it is
insecure. A manifold $M$ is secure if every pair of points in $M$
is secure. Otherwise $M$ is insecure. If $M$ is secure, and there
is a uniform upper bound on the cardinality of minimal blocking
sets, then $M$ is uniformly secure. On the other hand, $M$ is {\em
totally insecure} if all pairs $x,y\in M$ are insecure.

There are relationships between the (in)security of a compact
riemannian manifold and its topological entropy, fundamental
group, flatness of the metric, etc. For instance, the uniform
security of $M$ implies the vanishing of its topological entropy
and the quasi-nilpotence of $\pi_1(M)$ \cite{BG06}. If, in
addition, $M$ has no conjugate points, then it is flat
\cite{BG06}. The current conjecture is that a compact riemannian
manifold is uniformly secure iff it is flat \cite{BG06,LS}.

Set $m_t(x,y)=|\Gamma_t(x,y)|$ and $n_t(x,y)=|G_t(x,y)|$. These
are the {\em counting functions} for geodesics in $M$. Burns and
Gutkin \cite{BG06} related the security of $M$ with the growth of
counting functions as $t\to\infty$. The goal of the present paper
is to establish relationships between the growth of functions
$n_t(x,y)$ and the asymptotics of $s_t(x,y)$. This extends the
approach of \cite{BG06}. To see this, observe that i) a pair
$x,y\in M$ is secure iff $s_t(x,y)$ is a bounded function on
$\R_+$; ii) a manifold $M$ is uniformly secure iff the the
function $s(t)$ is bounded.

We will relate the topological entropy, volume entropy, and growth
rate of $\pi_1(M)$ with the asymptotics of functions $s_t(x,y)$
and $s(t)$. Before stating  our results, we need to say a few
words about infinite blocking costs. For almost all $x,y\in M$ we
have $n_t(x,y)<\infty$ \cite{BeBo}, and hence $s_t(x,y)<\infty$ as
well. Examples show that for some $x,y\in M$ and $t>t_0(x,y)$ we
may have $m_t(x,y)=\infty$ and $s_t(x,y)=\infty$,  implying
$s(t)=\infty$. In view of this,  we will often make provisos for
the possibility of infinite blocking costs.
The following proposition is a combination of
Theorem~\ref{fund_gr_thm} and Theorem~\ref{entro_cor} in
section~\ref{main} below.

\begin{theorem}   \label{intro_thm}
Let $M$ be a compact riemannian manifold.

\noindent i) If $\pi_1(M)$ grows exponentially, then the blocking
cost is either infinite or grows (at least) exponentially.

\noindent ii) Let $e>0$ be the topological entropy or the volume
entropy of $M$. Then the blocking cost is either infinite or grows
exponentially, with the  rate greater than or equal to $e/2$.
\end{theorem}

\medskip



The plan of the paper is as follows. In section~\ref{back&prev} we
expose the background material; we also sketch a proof of
Proposition~\ref{main_prop} which is our key technical result. Our
proof explains the relationship between the counting and blocking
of geodesics.
In section~\ref{mainsect} we expose several auxiliary propositions
that we will need in section~\ref{main}.
In section~\ref{main} we present our main results.


\noindent{\bf Notes and acknowledgements.} I thank Keith Burns for
comments on a draft of the paper.

\section{Background and Preliminaries} \label{back&prev}
It is convenient to partition the exposition into several
subsections.
\subsection{Counting geodesics between points, topological entropy, and volume entropy} \label{SSn_T}
Let $M$ be a compact riemannian manifold, let $d\mu$ be the
riemannian measure, and let $h=h(M)$ be the {\em topological
entropy} of $M$. Let $\tM$ be the universal cover of $M$. For
$\tx\in\tM$ let $B(\widetilde x,t)$ be the ball of radius $t$ in
$\tM$ around $\tx$. The exponential growth rate for
$t\mapsto\hbox{\rm Vol}\, B(\tx,t)$ does not depend on $\tx$; this
is the {\em volume entropy} $\la=\la(M)$. The two entropies are
related by $\la\le h$, and we have
\begin{equation}   \label{mane_eq}
h=\lim_{t\to\infty}\frac{1}{t}\log\int_{M\times
M}n_t(x,y)\,d\mu(x)d\mu(y).
\end{equation}
If $M$ has no {\em conjugate points}, then for any $x,y \in M$ we
have
\begin{equation}        \label{rate_eq}
h=\lim_{t \to \infty} \frac1t \log n_t(x,y).
\end{equation}
See  \cite{Man} and \cite{Mane} for this material.
Equation~\eqref{rate_eq} fails, in general, if $M$ has conjugate
points \cite{BP96,BP97}. From the obvious inequality
$$
s_t(x,y)\le m_t(x,y) \le n_t(x,y)
$$
and results of \cite{BeBo}, we have $s_t(x,y)<\infty$ for almost
all (resp. all) $x,y\in M$ (resp. if $M$ has no conjugate points).
The manifold is called secure (resp. uniformly secure, resp.
totally insecure) if $t\mapsto s_t(x,y)$ is bounded for all
$x,y\in M$ (resp. $t\mapsto s(\cdot)$ is bounded, resp. $t\mapsto
s_t(x,y)$ is unbounded for all $x,y\in M$). The following
proposition combines results of \cite{BG06,LS,GS06}.


\vspace{2mm}
\begin{theorem}   \label{combin_thm}
Let $M$ be a compact riemannian manifold.

\medskip

\noindent 1. If $M$ is uniformly secure, then it has zero
topological entropy and virtually nilpotent fundamental group. If,
in addition, $M$ has no conjugate points then it is flat.

\noindent 2. If $M$ has no conjugate points and positive
topological entropy, then it is totally insecure.

\noindent 3. Let  $M$ be a locally symmetric space. Then $M$ is
secure iff it is uniformly secure iff it is of euclidean type.
\end{theorem}

\medskip

These results suggest the following conjecture \cite{BG06,LS}.

\begin{conj}  \label{flat_conj}
A compact riemannian manifold is secure iff it is uniformly secure
iff it is flat.
\end{conj}

\subsection{Estimating the number of geodesics between points via the blocking cost}    \label{SSbl_cost}
For the benefit of the reader, we will sketch a proof of
Proposition~\ref{main_prop} which estimates the counting function
for geodesics via the blocking cost.

Let $B(x,y;t)$ be a minimal blocking set for $\Ga_t(x,y)$. Every
$\ga\in\Ga_t(x,y)$ is a concatenation:
$\ga=\{\al\in\Ga(x,z)\}\cup\{\be\in\Ga(z,y)\}$, where $z\in
B(x,y;t)$. Carefully choosing the point $z\in B(x,y;t)$, we obtain
the bound $m_t(x,y)\le\sum_{(p,q)\in\pp_1}m_{t/2}(p,q)$; the set
$\pp_1\in M\times M$ consists of pairs $(x,z),(z,y)$, where $z\in
B(x,y;t)$. Hence $|\pp_1|\le 2s(t)$.

Iterating this argument, we obtain a sequence of finite sets
$\pp_k\in M\times M$ and bounds
$m_t(x,y)\le\sum_{(p,q)\in\pp_k}m_{t/2^k}(p,q)$, where
$$
|\pp_k|\le 2^k\,s(t)\times\cdots\times s(\frac{t}{2^{k-1}}).
$$
The inductive process stops when $t/2^k$ gets smaller than the
injectivity radius of $M$; we then have $m_{t/2^k}(p,q)\le 1$ for
any pair $p,q\in M$. Let $\ka(t)$ be the smallest $k\in\N$
satisfying this inequality, and set $S(t)=s(t)\times\cdots\times
s(\frac{t}{2^{\ka(t)-1}})$. Then we have the bound
$$
m_t(x,y)\le 2^{\ka(t)}\,S(t).
$$
Since $\ka(t)\le \const\,\log_2t$, and $n_t(x,y)\le\const\,
t^2m_t(x,y)$, by Lemma~3.1 in \cite{BG06}, we have obtained a
desired bound.

\medskip

\subsection{Amplifications} \label{SS_rel}
The framework of security makes sense for any space with rich sets
of distinguished curves joining arbitrary pairs of points in the
space. In particular, it is meaningful for riemannian manifolds
with boundaries and corners. Planar billiard tables yield
elementary examples of this setting; billiard orbits play the role
of riemannian geodesics.

A polygon $P\subset\RR$ is secure if all billiard orbits
connecting an arbitrary pair of points in $P$ can be blocked by a
finite set. Which polygons are secure and which are insecure? A
complete answer is unknown but there are a few partial results
\cite{Gut05,Gut06}. For instance, the regular $n$-gon $R_n$ is
secure iff $n=3,4,6$ \cite{Gut05}. The counting functions
$n_t(x,y)$ are subexponential for polygons \cite{GH}; it is widely
believed that they are, in fact, polynomial \cite{Gut03}. Thus,
the results on insecurity of polygons in \cite{Gut05,Gut06} are
obtained using entirely different considerations.

P. Herreros studied the security for pairs of points in a $C^1$
riemannian surface $M$ homeomorphic to the $2$-sphere \cite{PH}.
Herreros found a large set of secure pairs of points whose
blocking sets have unexpected properties. This phenomenon has
applications to the security in riemannian products $M\times N$.

\section{Blocking cost and growth of joining geodesics}   \label{mainsect}
Let $M$ be a compact riemannian manifold. For $x,y\in M$ and $0<t$
let $s_t(x,y)$ be the minimal number of points needed to block all
$\ga\in\Ga_t(x,y)$, and set
\begin{equation}  \label{sec_cost_eq}
s(t)=\sup_{x,y\in M}s_t(x,y)\le\infty.
\end{equation}
Then $s_t(x,y)$ is the {\em blocking threshold} for the triple
$(x,y;t)$, and $s(t)$ is the {\em blocking cost}  function for
$M$.  If $s(t)=\infty$ for some $t$, we will say that the {\em
blocking cost is infinite}.\footnote{If the blocking threshold for
a particular $(x,y;t)$ is infinite, then the blocking cost is
infinite. Although it is apriori possible that $s_t(x,y)<\infty$
for all $(x,y;t)$ but the blocking cost is infinite, no examples
are known to the author.}


\vspace{3mm}

We will now introduce an operation on functions in $\R_+$;
although it is defined in a greater generality, we will restrict
our attention to positive, non-decreasing functions. The operation
depends on a positive parameter, $\de$, whose value will be set
later on. For the moment, $\de$ is arbitrary, and we suppress it
from notation.

For $t\in\R_+$ let $\ka(t)\in\N$ be the smallest $k$ such that
$\frac{t}{2^k}<\de$. Equivalently, $\ka(t)=0$ if $t<\de$ and
$\ka(t)=1+\lfloor\log_2t-\log_2\de\rfloor$ if $\de\le t$. We will
denote functions on $\R_+$ before they are fed into our operation
by $f,g$ etc; we use $F,G$ etc for the respective functions
produced by the operation.

Let $f$ be a function on $\R_+$. We set
\begin{equation}   \label{mult_eq}
F(t)=\prod_{0\le k \le \ka(t)-1}f(\frac{t}{2^k}).
\end{equation}
\begin{prop}   \label{main_prop}
Let $M$ be a compact riemannian manifold with finite blocking cost
$s(\cdot)$, and let $\de$ be its injectivity radius. Let
$S(\cdot)$ be the function associated with $s(\cdot)$ by
equation~\eqref{mult_eq}. Then for any $x,y\in M$ and $0<t$ we
have
$$
n_t(x,y) \le  \frac{t^3}{2\de^3} \,S(t).
$$
\begin{proof}
For $p,q\in M$ and $0<t$ denote by $B_t(p,q)\subset
M\setminus\{p,q\}$ a minimal blocking set for $\Ga_t(p,q)$.

For $t\in\R_+$ and any pair $x,y\in M$ we will define a sequence
$\pp_k$ of finite subsets of $M\times M$, where $0\le k \le
\ka(t)$.



If $t<\de$ then $\ka(t)=0$ and $\pp_0=\{(x,y)\}$. This defines the
sequence $\pp_k$ in this case. Let $\de\le t$. Set
$\pp_1=\{(x,z):z\in B_t(x,y)\}\cup\{(z,y):z\in B_t(x,y)\}$.


Suppose that $\pp_1,\dots,\pp_k$ have been defined. If
$\frac{t}{2^k}<\de$ then $k=\ka(t)$, and we have terminated the
sequence of sets. Otherwise we define $\pp_{k+1}\subset M\times M$
as the set comprised by pairs $(p,z),(z,q)$ where $(p,q)\in\pp_k$
and $z\in B_{t/2^k}(p,q)$.


By the argument of Lemma~3.2 in \cite{BG06}, for any $0\le k \le
\ka(t)$ we have
\begin{equation}   \label{reduc_eq}
m_t(x,y) \leq \sum_{(p,q)\in \pp_k}m_{t/2^k}(p,q).
\end{equation}
Set $\ppp=\ppp(x,y;t)=\pp_{\ka(t)}(x,y)$.
Applying equation~\eqref{reduc_eq} to $k=\ka(t)$ and taking into
account that $m_s(p,q)\le 1$ if $s<\de$, we obtain
\begin{equation}   \label{bound_eq1}
m_t(x,y) \leq |\ppp|.
\end{equation}
%

Minimal blocking sets are not unique, in general. Thus the set
$\ppp\subset M\times M$ is not uniquely defined, but we will use a
bound on $|\ppp|$.


We have $|\pp_0|=1$; for $1 \le k \le \ka(t)$ we have, by
construction
\begin{equation}   \label{cardin_eq}
|\pp_k(x,y;t)|\le 2s(t)\times 2s(\frac t2)\times\cdots\times
2s(\frac{t}{2^{k-1}})=2^kS_k(t).
\end{equation}
Since $S_{\ka(t)}(\cdot)=S(\cdot)$, and $\ka(t)\le
1+\log_2t-\log_2\de$, equation~\eqref{cardin_eq} yields
\begin{equation}   \label{bound_eq2}
m_t(x,y) \leq \frac{2t}{\de}\,S(t).
\end{equation}

By Lemma~3.1 in \cite{BG06}
\begin{equation}   \label{bu_gu_eq}
m_t(x,y) \le n_t(x,y)\le \frac{t^2}{4\de^2}\, m_t(x,y).
\end{equation}
Combining equations~\eqref{bound_eq2},~\eqref{bu_gu_eq}, we obtain
the claim.
\end{proof}
\end{prop}



For functions on $\R_+$, we will use the standard notation
$f=O(g),f=o(g),f\sim g$. The latter means  that there are
constants $0<c\le C<\infty$ such that
$$
c \, f(t) \le g(t) \le C\, f(t).
$$

We define the {\em rate of exponential growth} of a function by
$$
r=r(f)=\limsup_{t\to\infty}\frac{\log\,f(t)}{t}.
$$
Equivalently, $ r(f)=\inf\{a\in\R_+:\,f=O(e^{at})\}$. If
$0<r<\infty$ (resp. $r=0,\,r=\infty$) then $f$ grows exponentially
(resp. subexponentially, super-exponentially).


If $f=O(t^r)$ for some $0\le r$, we say that $f$ {\em grows (at
most) polynomially}. If $f\ne O(t^r)$ for any $r\in\R_+$, we say
that $f$ {\em is super-polynomial}.





%
\begin{lem}   \label{equiv_lem}
Let $f,g$ be functions on $\R_+$, and let $F,G$ be the functions
obtained from them via equation~\eqref{mult_eq}.

\noindent i) If $g=O(f)$ then there exists $\al\in\R$ such that
$G=O(t^{\al}F)$.

\noindent ii) If $g=o(f)$ then $G=O(t^{\be}F)$ for any $\be\in\R$.

\noindent iii) If $f\sim g$ then there exist $\al,\be\in\R$  such
that $G=O(t^{\al}F),\,F=O(t^{\be}G)$.
\begin{proof}
We use equation~\eqref{mult_eq} and the inequality
\begin{equation}   \label{estim_eq}
\log_2t-\log_2\de \le \ka(t) \le \log_2t-\log_2\de +1.
\end{equation}
Directly from the definition we see that $g=f_1f_2$ implies
$G=F_1F_2$ and that $f\le g$ implies $F\le G$.


Claim i) follows from preceding observations and from the
calculation of correspondence $f\to F$ when $f$ is a constant
function. We leave the latter to the reader.


Let now $g=o(f)$. Then for any $0<\ve$ we have the representation
$f=\vap_{\ve}\cdot g$ where $\vap_{\ve}(t)<\ve$ when $t_{\ve}<t$.
Let $0<c$ be the maximum of $\vap_{\ve}$ on $[0,t_{\ve}]$. Set
$\psi_{\ve}(t)=c$ if $t\le t_{\ve}$ and $\psi_{\ve}(t)=\ve$ if
$t_{\ve} <t$. Then $f\le\psi_{\ve}\cdot g$. Applying preceding
remarks and directly calculating the correspondence
$\psi_{\ve}\to\Psi_{\ve}$, we obtain ii).


Finally, iii) is a direct consequence of i), since $f\sim g$ means
that $f=O(g),g=O(f)$.
\end{proof}
\end{lem}
\medskip


We will use Proposition~\ref{main_prop} to relate the growth of
counting functions $n_t(x,y)$, as $t\to\infty$, and the
asymptotics at infinity of the blocking cost. We will need a
technical lemma.

\begin{lem}   \label{bound_lem}
Let $f$ be a function on $\R_+$, and let $F$ be the function
associated with it by equation~\eqref{mult_eq}.
Then the following statements hold.\\

\noindent i) If $f=O(e^{at})$ then for any $0<\ve$ we have
$F=O(e^{(2a+\ve)t})$.

\noindent ii) Let  $f=O(t^r)$ where $0\le r$. Then there exists
$\al\in\R$ such that $F=O(t^{\al\log_2t})$.
\begin{proof}
We have $f(t)< c\cdot e^{at}$ for some $0<c$. A direct calculation
from equation~\eqref{mult_eq} yields $e^{at}\to e^{2at}$.
Computing $F$ when $f=\const$, and using the argument of
Lemma~\ref{equiv_lem}, we obtain $F=O(t^{\al}e^{2at})$ for some
$\al\in\R$. This implies i).


We compute $F$ when $f(t)=t$. Directly from
equation~\eqref{mult_eq}, we obtain
$F(t)=t^{\ka(t)}2^{-\ka(t)(\ka(t)+1)/2}$. Estimating $\ka(t)$ via
equation~\eqref{estim_eq}, we obtain $F=O(t^{a+\frac 12
\log_2t})$. Whatever is the value of $a\in\R$, we have
$F=O(t^{(\ve+\frac 12)\log_2t})$ for any $0<\ve$.


From $f=O(t^r)$ we have $f(t)<\const\, t^n$ for some $n\in\N$. By
preceding remarks, this implies $F=O(t^a\cdot t^{n(\ve+\frac
12)\log_2t})$ for some $a\in\R$. Setting $\al=(2n+1)/2$, we obtain
ii).
\end{proof}
\end{lem}


Let $f,g$ be positive functions on $\R_+$. If $f=O(g)$ (resp.
$f=o(g)$), we say that $f$ grows not faster (resp. slower) than
$g$.

%
\begin{lem}   \label{super_pol_prop}
Let $M$ be a compact riemannian manifold. Suppose that for any
$0<c$ there exist $x,y\in M$ such that $n_t(x,y)$ grows faster
than $t^{c\log t}$. Then the blocking cost is either infinite or
super-polynomial.
\begin{proof}
Assume the opposite, i. e. that $s=O(t^r)$ for some $0<r$. Then,
by Lemma~\ref{bound_lem}, $S=O(t^{a\log t})$ for some $0<a$.
Proposition~\ref{main_prop} yields that $n_t(x,y)=O(t^{c\log t})$
for any $c$ greater than $a$,  and arbitrary $x,y\in M$. This
contradicts to the assumption.
\end{proof}
\end{lem}





\section{Blocking thresholds and geometry of a manifold}  \label{main}
Let $G$ be a finitely generated group, and let $W(G,S,n)\subset G$
be the set of elements of length at most $n$ with respect to a
finite generating set $S$. Denote by $r_S(G)$ the growth rate of
the function $w_S(n)=|W(G,S,n)|$. In general, $r_S(G)$ depends on
the choice of $S$. We say that $r_S(G)$ is the {\em growth rate of
$G$ with respect to the set $S$ of generators}. Let $S',S''$ be
two finite sets of generators, and let $r_{S'}(G),r_{S''}(G)$ be
the corresponding rates.  Then $0<r_{S'}(G)$ iff $0<r_{S''}(G)$.
Thus, we can speak of {\em groups with exponential growth} (resp.
{\em groups with sub-exponential growth}) without specifying the
growth rate.


The proposition below shows that the blocking cost function
controls the topological entropy of a manifold, as well as the
growth rate of the fundamental group.

\begin{prop}     \label{expo_thm}
Let $M$ be a compact riemannian manifold, with a finite blocking
cost $s(\cdot)$.

\noindent i) Let $0<\si$. If $s=O(e^{\si t})$ then $h(M)\le 2\si$.

\noindent ii) If the function $s(\cdot)$ grows subexponentially,
then the group $\pi_1(M)$ grows subexponentially as well.
\begin{proof}
i) If $f=O(e^{\si t})$, then, by Lemma~\ref{bound_lem},
$F=O(e^{(2\si+\ve) t})$ for any $0<\ve$. By Lemma~\ref{equiv_lem}
and Proposition~\ref{main_prop},  $n_t(x,y)=O(e^{(2\si+\ve) t})$
for all $x,y\in M$. By  Ma\~ne's formula equation~\eqref{mane_eq},
we have $h(M)\le 2\si+\ve$. Since $\ve$ is arbitrary, we obtain
the claim.

%
%


\noindent ii) By \cite{Man1,Miln}, $\pi_1(M)$ grows
subexponentially iff $\la(M)=0$. By i) and Manning's inequality
\cite{Man}, if $s(\cdot)$ is subexponential then $\la(M)=0$.
%
%
\end{proof}
\end{prop}

We will now turn to the main results of this paper.

\begin{thm}     \label{fund_gr_thm}
Let $M$ be a compact manifold. If the fundamental group of $M$
grows exponentially then the blocking cost for any riemannian
metric on $M$ is either infinite or grows at least  exponentially.
\begin{proof}
If the blocking cost function is subexponential, then, by
Proposition~\ref{expo_thm}, $\pi_1(M)$ grows subexponentially.
\end{proof}
\end{thm}


\begin{thm}     \label{entro_cor}
Let $M$ be a compact riemannian manifold. Let $\la=\la(M)$ and
$h=h(M)$ be its volume entropy and the topological entropy
respectively.

\noindent i) If $0<h$ then the blocking cost is either infinite or
grows exponentially, with the rate at least $h/2$.

\noindent ii) If $0<\la$ then the blocking cost is either infinite
or grows exponentially, with the rate at least $\la/2$.
\begin{proof}
Let $r$ be the exponential growth rate of $s(\cdot)$. By
Proposition~\ref{expo_thm}, $r\ge h/2$. This proves i). Claim ii)
follows from i) by Manning's inequality.


\end{proof}
\end{thm}


\begin{thm}     \label{vol_gr_cor}
Let $M$ be a compact riemannian manifold. Suppose that there is
$\tx\in\tM$ such that $\hbox{\rm Vol}\,B(\tx,t)$ grows faster than
any $t^{c\log t}$.
Then the blocking cost for $M$ is either infinite or grows
super-polynomially.
\begin{proof}
Assume the opposite, i. e., that $s=O(t^r)$ for some $0<r$. By
Lemma~\ref{bound_lem} and Proposition~\ref{main_prop}, there
exists $c\in\R_+$ such that $n_t(x,y)=O(t^{c\log t})$ for all
$x,y\in M$. By the proof of Proposition~\ref{expo_thm}, this
implies $\hbox{\rm Vol}\, B(\tx,t)=O(t^{c\log t})$ for all
$\tx\in\tM$, contrary to the assumption.
\end{proof}
\end{thm}
%

%


Preceding statements concern the blocking cost function, as
opposed to individual blocking thresholds $s_t(x,y)$. We will now
obtain estimates for those.

\begin{cor}     \label{top_ent_cor}
Let $M$ be a compact riemannian manifold. Let $h$ be its
topological entropy. If $0<h$, then for any $C>1$, arbitrarily
small $\ve>0$, and arbitrarily large $t$ there exist points
$x,y\in M$ such that
$$
Ce^{(\frac{h}{2}-\ve)t} \le s_t(x,y).
$$
\begin{proof}
Assume the opposite, i. e. that there is $C>1$, $\ve>0$ and
$\tau\in\R_+$ such that for all $x,y\in M$ and all $t$ greater
than $\tau$ we have the bound $s_t(x,y) <
Ce^{(\frac{h}{2}-\ve)t}$. Thus, $s=O(e^{(\frac{h}{2}-\ve)t})$.
Hence, the exponential growth rate of the blocking cost function
is strictly less than $h/2$. By Theorem~\ref{entro_cor}, this is
impossible.
\end{proof}
\end{cor}


\begin{cor}     \label{vol_ent_cor}
Let $M$ be a compact riemannian manifold. Let $\la$ be its volume
entropy. If $0<\la$, then for any $C>1$, arbitrarily small
$\ve>0$, and arbitrarily large $t$ there exist points $x,y\in M$
such that
$$
Ce^{(\frac{\la}{2}-\ve)t} \le s_t(x,y).
$$
\begin{proof}
The claim follows from Corollary~\ref{top_ent_cor} and the
inequality $\la\le h$ \cite{Man}.
\end{proof}
\end{cor}


\begin{cor}     \label{fund_gr_cor}
Let $M$ be a compact manifold whose fundamental group grows
exponentially. We endow $M$ with a riemannian metric. For $x,y\in
M$ and $t\in\R_+$ let $s_t(x,y)$ be the blocking threshold with
respect to this metric.

Then for any $r>1$ and arbitrarily large $t$ there exist points
$x,y\in M$ such that
$$
t^{r} < s_t(x,y).
$$
\begin{proof}
Suppose that the claim is false. Then there exist $r>1$ and
$\tau\in\R_+$ such that $s_t(x,y)\le t^{r}$ for all $\tau<t$ and
$x,y\in M$. Then the blocking cost $s(\cdot)$ is polynomial. This
contradicts Theorem~\ref{fund_gr_thm}.
\end{proof}
\end{cor}


\begin{rem}     \label{fund_gr_rem}
{\em Corollary~\ref{fund_gr_cor} says that if the fundamental
group of a manifold grows exponentially, then its blocking
thresholds are super-polynomial. By our methods we can obtain
other statements of that nature. For instance, if $\pi_1(M)$ grows
faster than any $n^{a\log n}$, then the blocking thresholds for
any riemannian metric on $M$ are super-polynomial. This is proved
by combining Lemma~\ref{bound_lem} and Lemma~\ref{super_pol_prop}.
An analogous statement holds if we replace the growth of
$\pi_1(M)$ by that of the volume of balls in $\tM$. We leave
details to the reader.}
\end{rem}


In order to illustrate the preceding material, we will now derive
some of the results of \cite{BG06} and \cite{LS}.\footnote{Compare
theorem~\ref{old_thm} with theorem~4.3 in \cite{BG06}.}
\begin{thm}     \label{old_thm}
Let $M$ be a uniformly secure compact riemannian manifold. Then
the topological entropy for the geodesic flow on $M$ vanishes. The
fundamental group of $M$ is virtually nilpotent.
\begin{proof}
In the present terminology, the blocking cost function for $M$ is
bounded. Hence, by Proposition~\ref{expo_thm}, $h(M)\le 2a$ for
any positive $a$. This proves the first claim.


Let $f$ be a bounded function on $\R_+$, and let $F$ be the
function corresponding to $f$ via equation~\eqref{mult_eq}. As we
have seen in the proof of Lemma~\ref{equiv_lem}, if $f$ is bounded
then $F$ grows polynomially. Applying this remark to the blocking
cost function of $M$ and using Proposition~\ref{main_prop}, we
obtain a uniform bound $n_t(x,y)< \const\,t^r$. By the proof of
Proposition~\ref{expo_thm}, the group $\pi_1(M)$ has polynomial
growth. Our second claim now follows from a theorem of Gromov
\cite{Grom}.
\end{proof}
\end{thm}


\begin{rem}   \label{no_conj_pts_rem}
{\em If the manifold in Theorem~\ref{old_thm} has no conjugate
points, then the conclusion is much stronger. Namely, by a theorem
of Lebedeva, a compact riemannian manifold with no conjugate
points and a quasi-nilpotent fundamental group is flat \cite{Leb}.
See theorem~4.3 in \cite{BG06}. This is one of the results
supporting Conjecture~\ref{flat_conj}.}
\end{rem}


The following proposition strengthens Theorem~4.5 in
\cite{BG06}.\footnote{Theorem~4.5 has been independently obtained
by Lafont and Schmidt \cite{LS}.}

\begin{thm}   \label{nocp_thm}
Let $M$ be a compact riemannian manifold with no conjugate points.
Let $h>0$ be its topological entropy. Then for any $x,y\in M$ the
blocking threshold $s_t(x,y)$ grows exponentially. Its exponential
growth rate is at least $h/2$.
\end{thm}
\begin{proof}
Let $x,y\in M$ be arbitrary, and let $t\in\R_+$.
The beginning of the proof of Proposition~\ref{main_prop} yields
\begin{equation}             \label{half_time_eq1}
 m_t(x,y) \le \sum_{(p,q)\in\pp_1}  m_{t/2}(p,q),
\end{equation}
where $|\pp_1|=2s_t(x,y)$. Let $\de>0$ be the injectivity radius
of $M$. By equations~\eqref{rate_eq},~\eqref{bu_gu_eq}, we
have\footnote{Here we use that $M$ has no conjugate points.
Convergence in equation~\eqref{con_entr_eq} is uniform
\cite{Mane}.}
\begin{equation} \label{con_entr_eq}
\lim_{t\to \infty} \frac1t \log m_t(p,q) = h.
\end{equation}
%


Let $\ve>0$ be arbitrary. Combining
equations~\eqref{half_time_eq1},~\eqref{con_entr_eq}, we obtain
for $t>t(\ve)$,
\begin{equation} \label{get_there_eq}
e^{(h-\ve)t}\le m_t(x,y) \le 2s_t(x,y)e^{\frac12(h+\ve)t}.
\end{equation}
%
Equation~\eqref{get_there_eq} yields $ e^{(\frac
h2-\frac{3\ve}{2})t} \le 2s_t(x,y)$; letting $\ve\to 0$, we obtain
the claim.
\end{proof}



\medskip

Applying Theorem~\ref{nocp_thm} to manifolds with nonpositive
curvatures, we obtain the following proposition. It strengthens
Corollary~4.7 in \cite{BG06}.

\begin{cor}     \label{neg_curv_cor}
Let $M$ be a compact riemannian manifold of nonpositive curvature.
Then the following dichotomy holds:

\noindent i) The manifold is uniformly secure. Its security
threshold is bounded above in terms of the dimension of $M$;

\noindent ii) The blocking thresholds $s_t(x,y)$ grow
exponentially. Their growth rates are greater than or equal to
half the topological entropy of $M$.
\end{cor}
\begin{proof}
By \cite{Goodw,Pes}, we have the dichotomy: i) $M$ is flat; ii)
$M$ has positive topological entropy. In case i), Proposition 2 in
\cite{GS06} yields the claim. In case ii), the claim follows from
Theorem~\ref{nocp_thm}.
\end{proof}


\begin{thebibliography}{99}
%


\bibitem{BeBo} M. Berger and R. Bott, {\em Sur les vari\'et\'es \`a courbure strictement positive},
Topology {\bf 1} (1962), 302 -- 311.





\bibitem{BG06} K. Burns and E. Gutkin, {\em Growth of the number of geodesics between points and insecurity
for Riemannian manifolds}, preprint math.DS/0701579 (2007).

\bibitem{BP96} K. Burns and G. Paternain, {\em On the growth
of the number of geodesics joining two points},
Pitman Res. Notes Math. {\bf 362} (1996), 7 -- 20.



\bibitem{BP97} K. Burns and G. Paternain, {\em
Counting geodesics on a Riemannian manifold and topological entropy of geodesic flows},
Erg. Theo. \& Dyn. Sys. {\bf 17} (1997), 1043 -- 1059.





















\bibitem{Grom} M. Gromov, {\em Groups of polynomial growth and expanding maps},
Inst. Hautes \'Etudes Sci. Publ. Math. {\bf 53} (1981), 53 -- 73.

\bibitem{Goodw} L.W. Goodwyn,  {\em Topological entropy bounds measure-theoretic entropy},
Proc. AMS  {\bf 23}  (1969), 679 -- 688.



\bibitem{Gut03} E. Gutkin, {\em Billiard dynamics: A survey with the emphasis
on open problems}, Reg. \& Chaot. Dyn. {\bf 8} (2003), 1 -- 13.





\bibitem{Gut05} E. Gutkin, {\em Blocking of billiard orbits and security for polygons and flat surfaces},
GAFA: Geom. \& Funct. Anal. {\bf 15} (2005), 83 -- 105.

\bibitem{Gut06} E. Gutkin, {\em Insecure configurations in lattice translation surfaces,
with applications to polygonal billiards}, Discr. Cont. Dyn. Sys. {\bf A 16} (2006), 367 -- 382.

\bibitem{GH} E. Gutkin and N. Haydn, {\em Topological entropy
of polygon exchange transformations and polygonal billiards.},
Erg. Theo. \& Dyn. Sys. {\bf 17} (1997), 849 -- 867.





\bibitem{GS06} E. Gutkin and V. Schr\"oder, {\em Connecting geodesics and
security of configurations in compact
locally symmetric spaces}, Geometriae Dedicata {\bf 118} (2006), 185 -- 208.






\bibitem{PH} P. Herreros, {\em Blocking: New examples and properties of
products}, preprint    arXiv:0707.0456 (2007).







\bibitem{LS} J.-F. Lafont and B. Schmidt, {\em Blocking light in compact Riemannian manifolds},
Geometry \& Topology {\bf 11} (2007), 867 -- 887.

\bibitem{Leb}
N.D. Lebedeva, {\em  On spaces of polynomial growth with no conjugate points},
St. Petersburg Math. J.  {\bf 16}  (2005), 341 -- 348.



\bibitem{Mane} R. Ma\~n\'e, {\em On the topological entropy of geodesic flows},
J. Diff. Geom.  {\bf 45} (1997), 74 -- 93.

\bibitem{Man} A. Manning, {\em Topological entropy for geodesic flows}, Ann. Math.
{\bf 110}  (1979), 567 -- 573.

\bibitem{Man1} A. Manning, {\em Relating exponential growth in a
manifold and its fundamental group}, Proc. AMS {\bf 133}  (2004),
995 -- 997.

\bibitem{Miln} J. Milnor, {\em A note on curvature and fundamental group},
J. Diff. Geom.  {\bf 2}  (1968), 1 -- 7.









\bibitem{Pes} Ya. Pesin, {\em Formulas for the entropy of the geodesic flow on
a compact Riemannian manifold without conjugate points}, Math.
Notes {\bf 24} (1978), 796 -- 805.



























































\end{thebibliography}
\end{document}